\input ./style/arxiv-general.cfg
\documentclass[aap,MSNbibl,dvips]{arximspdf}
\makeatletter
   \@ifpackageloaded{graphicx}{}{\usepackage{graphicx}}
\makeatother

\doi{10.1214/14-AAP1066} 
\volume{25}
\issue{6}
\pubyear{2015}
\firstpage{3033}
\lastpage{3046}
\docsubty{FLA}

\makeatletter
\newcommand{\ds}{\displaystyle}
\let\geqslant\geq
\let\leqslant\leq
\newcommand{\rrvert}{\vert}
\newcommand{\rrVert}{\Vert}
\newcommand{\llvert}{\vert}
\newcommand{\llVert}{\Vert}

\newcommand{\eqref}[1]{(\ref{#1})}
\newtheorem{lemma}{Lemma}

\makeatother

\begin{document}
\begin{frontmatter}

\title{A diffusion process associated with  Fr\'echet means}
\runtitle{Diffusion associated with Fr\'echet means}

\begin{aug}
\author[A]{\fnms{Huiling}~\snm{Le}\corref{}\ead[label=e1]{huiling.le@nottingham.ac.uk}}
\runauthor{H. Le}
\affiliation{University of Nottingham}
\address[A]{School of Mathematical Sciences\\
University of Nottingham\\
Nottingham, NG7 2RD\\
United Kingdom\\
\printead{e1}}
\end{aug}

\received{\smonth{7} \syear{2013}}
\revised{\smonth{9} \syear{2014}}

\begin{abstract}
This paper\vspace*{1pt} studies rescaled images, under $\exp^{-1}_\mu$, of the
sample Fr\'echet means of i.i.d. random variables $\{X_k\vert
k\geqslant 1\}$ with Fr\'echet mean $\mu$ on a Riemannian manifold. We show that,
with appropriate scaling, these images converge weakly to a diffusion
process. Similar to the Euclidean case, this limiting diffusion is a
Brownian motion up to a linear transformation. However, in addition to
the covariance structure of $\exp^{-1}_\mu(X_1)$, this linear
transformation also depends on the global Riemannian structure of the manifold.
\end{abstract}

\end{frontmatter}

\section{Introduction}\label{sec1}

It has become increasingly common in various research areas for
statistical analysis to involve data that lies in non-Euclidean spaces.
One such an example is the statistical analysis of shape; cf. \cite{DM}
and \cite{KBCL}. Consequently, many statistical concepts and
techniques have been generalised and developed to adapt to such phenomena.

Fr\'echet means, as a generalisation of Euclidean means, of random
variables on a metric space have been widely used for statistical
analysis of non-Euclidean data. A point $\mu$ in a metric space
$\mathbf{M}$
with distance function $\rho$ is called a Fr\'echet mean of a random
variable $X$ on $\mathbf{M}$ if it satisfies
\[
\mu=\operatorname{arg}\min_{x\in\mathbf{M}}\mathrm{E}\bigl[
\rho(x,X)^2\bigr].
\]
Influenced by the structure of the underlying spaces, Fr\'echet means,
unlike their Euclidean counterparts, exhibit many challenging
probabilistic and statistical features. Various aspects of Fr\'echet
means have been studied for non-Euclidean spaces, including Riemannian
manifolds and certain stratified spaces. Among others, the strong law
of large numbers for Fr\'echet means on general metric spaces was
obtained in \cite{HZ}. The first use of Fr\'echet means to provide
nonparametric statistical inference, such as confidence regions and
two-sample tests for discriminating between two distributions, was
carried out in \cite{BP1} and \cite{BP} for both extrinsic and
intrinsic inference applied to manifolds. When $\mathbf{M}$ is a
Riemannian manifold with the distance function being that induced by
its Riemannian metric, the results on central limit theorems for Fr\'echet means can be found in \cite{BP}
and \cite{KL}. The results in
both papers imply that, since manifolds are locally homeomorphic to
Euclidean spaces, the limiting distributions for sample Fr\'echet means
on Riemannian manifolds are usually Gaussian, a phenomenon similar to
that for Euclidean means.

In the case of Euclidean space, the link between the sample means of
i.i.d. random vectors and random walks leads to the fact that the
rescaled sample means converge weakly to Brownian motion, possibly up
to a linear transformation associated with the covariance structure of
the random vectors. On the other hand, the authors of \cite{ADPY}
constructed a stochastic gradient algorithm from a given sequence of
i.i.d. random variables on a Riemannian manifold where, under certain
conditions, the random sequence resulting from the algorithm converges
almost surely to the Fr\'echet mean $\mu$ of the given random
variables. Moreover, it showed that, if one rescales the images, under
$\exp^{-1}_\mu$, of the random walks associated with the algorithm,
they converge weakly to an inhomogeneous diffusion process on the
tangent space of the manifold at $\mu$. The following questions are
raised from this paper: if one rescales the images, under $\exp^{-1}_\mu$, of the sample Fr\'echet means of the random variable, will they
converge weakly? If they do, do they converge to the same diffusion
process as the one given in \cite{ADPY}? If not, what is the limiting
diffusion process? This paper addresses these questions. We show that
the rescaled images of the sample Fr\'echet means of i.i.d. random
variables $\{X_k\vert k\geqslant1\}$ on a Riemannian manifold converge
weakly to a diffusion process which is a Brownian motion up to a linear
transformation. Moreover, in addition to the covariance structure of
$\exp^{-1}_\mu(X_1)$, this linear transformation also depends on the
global Riemannian structure of the manifold. For this we first, in the
next section, construct a sequence of simpler inhomogeneous Markov
processes, each of which is also a martingale, and consider the
behaviour of their weak convergence. In addition to their intrinsic
interest, the results in this section also form a basis for our
investigations of ``rescaled'' sample Fr\'echet means in the following
section. In particular, we relate the constructed sequence of processes
to the ``rescaled'' sample Fr\'echet means in such a way that the result
for the latter is a direct consequence of the former.

\section{An auxiliary weakly convergent sequence of Markov chains}\label{sec2}

Let $\mathbf{M}$ be a complete Riemannian manifold of dimension $d$ with
covariant derivative $D$ and Riemannian distance $\rho$, whose
sectional curvature is bounded below by $\kappa_0\leqslant0$ and above
by $\kappa_1\geqslant0$. For any $x\in\mathbf{M}$, we denote by
$\mathcal{C}_x$ the cut locus of $x$. Note that, for any fixed $x_0$, the
squared distance function $\rho(x_0,x)^2$ to $x_0$ is not $C^2$ on
$\mathcal{C}_{x_0}$.

For a fixed $y\in\mathbf{M}$, consider the vector field on $\mathbf{M}\setminus
\mathcal{C}_y$ defined, at $x\notin\mathcal{C}_y$, by $\exp
_x^{-1}(y)\in\tau_x({\mathbf{M}})$, where $\tau_x({\mathbf{M}})$
denotes the
tangent space of $\mathbf{M}$ at $x$, and then define the linear operator
$H_{x,y}$ on the tangent space $\tau_x({\mathbf{M}})$ by
%
\begin{equation} \label{eqn2}
H_{x,y}\dvtx v\mapsto- \bigl(D_v\exp^{-1}_x(y)
\bigr) (x).
\end{equation}
The operator $H_{x,y}$ so defined will play an important role in the
following study of the asymptotic behaviour of sample Fr\'echet means
on $\mathbf{M}$. Note first that $H_{x,y}$ is closely linked with
$\operatorname{Hess}(\frac{1}{2}\rho(x,y)^2 )$, the Hessian of the function $\frac{1}{2}\rho(x,y)^2$, as follows (cf. \cite{JJ}, page 145):
\[
\bigl\langle H_{x,y}(v),u\bigr\rangle=\operatorname{Hess}_x
\bigl(\tfrac{1}{2}\rho (x,y)^2 \bigr) (v,u),
\]
for any $x\notin\mathcal{C}_y$ and any tangent vectors $u,v\in\tau
_x({\mathbf{M}})$, and so the assumption on the bounds for the sectional
curvature of $\mathbf{M}$ implies that, for any unit tangent vector
$v\in\tau
_x({\mathbf{M}})$,
%
\begin{eqnarray}
&&\sqrt{\kappa_1}\rho(x,y)\cot\bigl(
\sqrt{\kappa_1}\rho(x,y)\bigr)
\nonumber
\\[-8pt]
\label{eqn0}\\[-8pt]
\nonumber
&& \qquad \leqslant\bigl\langle H_{x,y}(v),v\bigr\rangle\leqslant \sqrt{-
\kappa_0}\rho(x,y)\coth\bigl(\sqrt{-\kappa_0}\rho(x,y)
\bigr),
\end{eqnarray}
where we require that if $\kappa_1>0$, $\sqrt{\kappa_1}\rho
(x,y)<\pi/2$
for the first inequality to hold; cf. \cite{JJ}, page 203.

In contrast to Euclidean means, there is generally no closed form for
Fr\'echet means. On the other hand, the result of \cite{LB} implies
that the Euclidean random variable $\exp_\mu^{-1}(X)$ is almost surely
defined, where $\mu$ is a Fr\'echet mean of the random variable $X$ on
$\mathbf{M}$. Then, since
\[
\exp^{-1}_x(y)=-\tfrac{1}{2}\operatorname{grad}_1
\bigl(\rho(x,y)^2 \bigr),
\]
where grad$_1$ denotes the gradient operator acting on the first
argument of a function on ${\mathbf{M}}\times{\mathbf{M}}$ and since
\[
\operatorname{grad} \bigl(\mathrm{E}\bigl[\rho(x,X)^2\bigr] \bigr)
|_{x=\mu
}=\mathrm{E}\bigl[\operatorname{grad}_1 \bigl(
\rho(x,X)^2 \bigr) |_{x=\mu}\bigr]=0
\]
by the definition of Fr\'echet means, the Fr\'echet mean $\mu$
satisfies the condition that
%
\begin{equation}\label{eqn1}
\mathrm{E}\bigl[\exp_\mu^{-1}(X)\bigr]=0.
\end{equation}
Thus $\mu$ is linked to the Euclidean mean in the sense that the origin
of the tangent space of $\mathbf{M}$ at $\mu$, $\tau_\mu(\mathbf{M})$ is the
Euclidean mean of the Euclidean random variable $\exp^{-1}_\mu(X)$.

Let $\{X_k\vert k\geqslant1\}$ be a sequence of i.i.d. random variables
on $\mathbf{M}$, and for a fixed $x_0\in\mathbf{M}$, assume that
$\mathrm{E}[\rho (x_0,X_1)^2]<\infty$. This assumption ensures
the existence of Fr\'
echet means of $X_1$. For simplicity, in the following, we shall assume
that the Fr\'echet mean $\mu$ of $X_1$ is unique. However, we do not
require the support of the probability measure of $X_1$ to be contained
in any geodesic ball. Note that the result of \cite{LB} ensures that
$\mathrm{P}(X_1\in\mathcal{C}_\mu)=0$ under some mild
condition on $\mathbf{M}$.
We further assume that
%
\begin{eqnarray}
&& \mathrm{E}\bigl[\|H_{\mu,X_1}
\|^2\bigr]<\infty \quad\mbox{and}
\nonumber
\\[-8pt]
\label{eqn3}\\[-8pt]
\nonumber
&& \bigl(\mathrm{E}[H_{\mu,X_1}] \bigr)^{-1}, \mbox{the
inverse of the linear operator }\mathrm{E}[H_{\mu,X_1}], \mbox{exists}.
\end{eqnarray}
These two assumptions ensure that the linear operator $\mathrm{E}[H_{\mu,X_1}]$
is well defined and nonsingular.

For each fixed $n\geqslant1$, consider the time-inhomogeneous Markov
chain $\{V^n_k\vert k\geqslant0\}$ defined on the tangent space $\tau
_\mu
({\mathbf{M}})$ in terms of $\{X_k\vert k\geqslant1\}$ as
\begin{eqnarray*}
V^n_0&=&0,
\\
V^n_1&=&\frac{1}{\sqrt n} \bigl(\mathrm{E}[H_{\mu,X_1}]
\bigr)^{-1} \bigl(\exp ^{-1}_\mu(X_1)
\bigr),
\\
V^n_{k+1}&=&\frac{1}{\sqrt n} \bigl(\mathrm{E}[H_{\mu
,X_1}]
\bigr)^{-1} \bigl(\exp ^{-1}_\mu(X_{k+1})
\bigr)
\\
&&{}+ \biggl\{\frac{k+1}{k}I-\frac{1}{k} \bigl(\mathrm{E}[H_{\mu
,X_1}]
\bigr)^{-1}H_{\mu,X_{k+1}} \biggr\}\bigl(V^n_k
\bigr), \qquad k\geqslant1.
\end{eqnarray*}
One may check that $\{V^n_k\vert k\geqslant0\}$ is a martingale. We are
interested in the asymptotic behaviour of $\{V^n_{[nt]}\vert t\geqslant
0\}$ as $n$ tends to infinity. Firstly, for this, the following lemma
gives an upper bound for the sequence $\{V^n_{[\varepsilon_0n]}\vert n\geqslant1\}$ for $\varepsilon_0>0$.

\begin{lemma}\label{lem1}
Suppose that the assumptions \eqref{eqn3} hold. Then there is a
constant $c_0>0$ such that, for any $\varepsilon_0>0$ and $n_0>0$,
\[
\sup_{n\geqslant n_0}\mathrm{E}\bigl[\bigl|V^n_{[\varepsilon
_0n]}\bigr|^2
\bigr]\leqslant\alpha\mathrm{E}\bigl[\rho (\mu,X_1)^2
\bigr] \biggl\{ \frac{1}{n_0}+\varepsilon_0c_0 \biggr\},
\]
where $\alpha=\llVert  (\mathrm{E}[H_{\mu,X_1}]
)^{-1}\rrVert ^2$.
\end{lemma}

\begin{pf}
Write
\[
B=\mathrm{E}\bigl[H_{\mu,X_1}^\top \bigl(
\mathrm{E}[H_{\mu
,X_1}] \bigr)^{-\top} \bigl(
\mathrm{E}[H_{\mu,X_1}] \bigr)^{-1}H_{\mu,X_1}\bigr].
\]
Then, by the definition of $V^n_k$,
\begin{eqnarray*}
&&\mathrm{E}\bigl[\bigl|V^n_{k+1}\bigr|^2\vert
V^n_k\bigr]
\\
&& \qquad =\frac{1}{n} \mathrm{E}\bigl[\bigl|\bigl(\mathrm{E}[H_{\mu
,X_1}]
\bigr)^{-1}\exp^{-1}_\mu (X_1)\bigr|^2
\bigr]+ \biggl\langle V^n_k, \biggl\{\frac{k^2-1}{k^2}I+
\frac
{1}{k^2}B \biggr\}V^n_k \biggr\rangle
\\
&&\qquad\quad{}+\frac{2}{\sqrt{n}} \mathrm{E}\bigl[ \bigl\langle\bigl(
\mathrm{E}[H_{\mu,X_1}] \bigr)^{-1}\exp^{-1}_\mu(
X_{k+1}),G_{k+1}V^n_k \bigr\rangle\vert
V^n_k\bigr],
\end{eqnarray*}
where
\[
G_{k+1}=\frac{k+1}{k}I-\frac{1}{k} \bigl(
\mathrm{E}[H_{\mu
,X_1}] \bigr)^{-1}H_{\mu,X_{k+1}}.
\]
Under the given conditions, there is a constant $\beta\geqslant0$ such
that, for any $v\in\tau_\mu({\mathbf{M}})$, $\langle v,Bv\rangle
\leqslant
(\beta+1)|v|^2$. Thus, using the facts that $\mathrm{E}[\exp
^{-1}_\mu(X_1)]=0$
and that $X_{k+1}$ is independent of $V^n_k$, we have
\begin{eqnarray*}
&&\mathrm{E}\bigl[\bigl|V^n_{k+1}\bigr|^2\vert
V^n_k\bigr]
\\
&& \qquad \leqslant \frac{\alpha}{n}\mathrm{E}\bigl[\rho(\mu ,X_1)^2
\bigr]+ \biggl(1+\frac{\beta
}{k^2} \biggr)\bigl|V^n_k\bigr|^2
\\
&&\qquad\quad{}-\frac{2}{\sqrt{n}k} \bigl\langle\mathrm{E}\bigl[H_{\mu
,X_{k+1}}^\top
\bigl(\mathrm{E}[H_{\mu,X_1}] \bigr)^{-\top
} \bigl(
\mathrm{E}[H_{\mu,X_1}] \bigr)^{-1}\exp ^{-1}_\mu(
X_{k+1})\bigr],V^n_k \bigr\rangle.
\end{eqnarray*}
Noting that $\{V^n_k\vert k\geqslant0\}$ is a martingale with zero
expectation, the above implies that
\[
\mathrm{E}\bigl[\bigl|V^n_{k+1}\bigr|^2\bigr]
\leqslant\frac{\alpha
}{n}\mathrm{E}\bigl[\rho(\mu,X_1)^2
\bigr]+ \biggl(1+\frac{\beta}{k^2} \biggr)\mathrm{E}\bigl[\bigl|V^n_k\bigr|^2
\bigr].
\]
Hence, by induction, we have
\[
\mathrm{E}\bigl[\bigl|V^n_k\bigr|^2\bigr]
\leqslant\frac{\alpha}{n}\mathrm{E}\bigl[\rho(\mu,X_1)^2
\bigr] \Biggl\{1+\sum_{i=1}^k\prod
_{j=i}^k \biggl(1+ \frac{\beta}{j^2} \biggr) \Biggr
\}.
\]
Since
\[
\sum_{i=1}^k\prod
_{j=i}^k \biggl(1+ \frac{\beta}{j^2} \biggr)
\leqslant k\prod_{j=1}^k \biggl(1+
\frac{\beta}{j^2} \biggr),
\]
the above implies, in particular, that
\[
\mathrm{E}\bigl[\bigl|V^n_{[\varepsilon_0n]}\bigr|^2\bigr]
\leqslant\alpha \mathrm{E}\bigl[\rho(\mu,X_1)^2
\bigr] \Biggl\{ \frac{1}{n}+\varepsilon_0 \prod
_{j=1}^\infty \biggl(1+\frac{\beta
}{j^2} \biggr) \Biggr\}.
\]
The required result then follows from the fact that $\prod_{j=1}^\infty (1+\frac{\beta}{j^2} )<\infty$.
\end{pf}

The next lemma gives various bounds on the differences
$V^n_{k+1}-V^n_k$ for sufficiently large $n$ and $k$.

\begin{lemma}\label{lem2}
In addition to the assumptions in \eqref{eqn3}, assume that, for some
$\delta>0$, $E[\rho(\mu,X_1)^{2+\delta}]<\infty$. Then, for any
$\varepsilon_0>0$ and $r>0$, there are constants $c_1$, $c_2$, and $c_3$
depending on $\varepsilon_0$ and $r$ such that, when $n$ is sufficiently
large, for $k\geqslant\varepsilon_0n$ and for $v\in\tau_\mu({\mathbf
{M}})$ with
$|v|\leqslant r$:
\begin{longlist}[(iii)]
\item[(i)] for any $\varepsilon>0$
%
\begin{equation}
\label{eqn4}
\mathrm{P}\bigl(\bigl\llvert V^n_{k+1}-V^n_k
\bigr\rrvert >\varepsilon\vert V^n_k=v\bigr)\leqslant \cases{
\ds\frac{c_1}{\varepsilon^{2+\delta}}n^{-(1+\min\{
1,\delta
/2\})}, \vspace*{4pt}&\quad$\mbox{if } \varepsilon
\leqslant 1$,
\cr
\ds \frac{c_1}{\varepsilon
^2}n^{-(1+\min\{1,\delta/2\})},& \quad$\mbox{if }\varepsilon>1$;}
\end{equation}
\item[(ii)]
%
\begin{equation} \label{eqn5}
\bigl\llvert \mathrm{E}\bigl[ \bigl(V^n_{k+1}-V^n_k
\bigr)1_{\{\llvert V^n_{k+1}-V^n_k\rrvert >1\}}\vert V^n_{k}=v\bigr]\bigr\rrvert
\leqslant c_2n^{-(1+\min\{1/2,\delta/4\})};
\end{equation}
\item[(iii)]
%
\begin{eqnarray}
&&\bigl\llVert \mathrm{E}\bigl[
\bigl(V^n_{k+1}-V^n_k\bigr)
\bigl(V^n_{k+1}-V^n_k
\bigr)^\top 1_{\{ |V^n_{k+1}-V^n_k|>1\}}\vert V^n_k=v\bigr]
\bigr\rrVert
\nonumber
\\[-8pt]
\label{eqn6}\\[-8pt]
\nonumber
&&\qquad\leqslant c_3n^{-(1+\min\{1,\delta/(2+\delta)\})}.
\end{eqnarray}
\end{longlist}
\end{lemma}

\begin{pf}
For any $\varepsilon>0$, write $\varepsilon'=\varepsilon\|\mathrm{E}[H_{\mu,X_1}]\|$.
Then, by the definition of~$V^n_k$, we have
\begin{eqnarray*}
&&\mathrm{P}\bigl(\bigl\llvert V^n_{k+1}-V^n_k
\bigr\rrvert >\varepsilon\vert V^n_k=v\bigr)
\\
&& \qquad\leqslant\mathrm{P}\biggl(\biggl\llvert \exp^{-1}_\mu(X_1)+
\frac{\sqrt
{n}}{k} \bigl\{\mathrm{E}[H_{\mu,X_1}]-H_{\mu,X_1}
\bigr\} (v)\biggr\rrvert >\varepsilon'\sqrt{n}\biggr)
\\
&& \qquad\leqslant \mathrm{P}\biggl(\bigl\llvert \exp^{-1}_\mu(X_1)
\bigr\rrvert >\varepsilon'\frac{\sqrt n}{2}\mbox{ or }\bigl\llvert
\bigl\{\mathrm{E}[H_{\mu,X_1}]-H_{\mu,X_1} \bigr\} (v)\bigr
\rrvert >\varepsilon'\frac
{k}{2}\biggr)
\\
&& \qquad \leqslant \mathrm{P}\biggl(\bigl\llvert \exp^{-1}_\mu(X_1)
\bigr\rrvert >\varepsilon'\frac{\sqrt n}{2}\biggr)+\mathrm{P}
\biggl(\bigl\llVert \mathrm{E}[H_{\mu,X_1}]-H_{\mu,X_1}\bigr
\rrVert >\varepsilon'\frac{k}{2|v|}\biggr).
\end{eqnarray*}
Thus if $v\neq0$, it follows from Chebyshev's inequality that
\begin{eqnarray*}
&&\mathrm{P}\bigl(\bigl\llvert V^n_{k+1}-V^n_k
\bigr\rrvert >\varepsilon\vert V^n_k=v\bigr)
\\
&& \qquad \leqslant \mathrm{E}\bigl[\rho(\mu,X_1)^{2+\delta}\bigr]
\frac
{2^{2+\delta}}{(\varepsilon
'\sqrt{n})^{2+\delta}}+\operatorname{var}\bigl(\|H_{\mu,X_1}\|\bigr)\frac
{(2|v|)^2}{(\varepsilon'k)^2}
\\
&& \qquad \leqslant \mathrm{E}\bigl[\rho(\mu,X_1)^{2+\delta}\bigr]
\frac
{2^{2+\delta}}{(\varepsilon
'\sqrt{n})^{2+\delta}}+\mathrm{E}\bigl[\|H_{\mu,X_1}\|^2\bigr]
\frac
{(2|v|)^2}{(\varepsilon' k)^2},
\end{eqnarray*}
when $n$ is sufficiently large. Note that the assumption that $v\neq0$
implies that $k\geqslant1$. If $v=0$, a modified argument will show
that the above still holds for $k\geqslant1$. Hence, \eqref{eqn4} follows.

Similarly, using the definition of $V^n_k$, we have
\begin{eqnarray*}
&&\mathrm{E}\bigl[\bigl\llvert V^n_{k+1}-V^n_k
\bigr\rrvert ^2\vert V^n_k=v\bigr]
\\
&& \qquad \leqslant \frac{2}{n}\mathrm{E}\bigl[\bigl\llvert \bigl(
\mathrm{E}[H_{\mu,X_1}] \bigr)^{-1}\exp ^{-1}_\mu(X_1)
\bigr\rrvert ^2\bigr]+\frac{2}{k^2}\mathrm{E}\bigl[\bigl
\llvert \bigl(I- \bigl(\mathrm{E}[H_{\mu
,X_1}] \bigr)^{-1}H_{\mu,X_1}
\bigr)v\bigr\rrvert ^2\bigr].
\end{eqnarray*}
Thus, under the given conditions, result (i) also implies that, for
any $r>0$ and some constant $c_2$ depending on $\varepsilon_0$ and $r$,
we have
\begin{eqnarray*}
&&\bigl\llvert \mathrm{E}\bigl[ \bigl(V^n_{k+1}-V^n_k
\bigr)1_{\{\llvert V^n_{k+1}-V^n_k\rrvert >1\}}\vert V^n_{k}=v\bigr]\bigr\rrvert
\\
&& \qquad \leqslant  \mathrm{E}\bigl[\bigl\llvert V^n_{k+1}-V^n_k
\bigr\rrvert ^2\vert V^n_k=v
\bigr]^{1/2} \mathrm{P}\bigl(\bigl\llvert V^n_{k+1}-V^n_k
\bigr\rrvert >1\vert V^n_{k}=v\bigr)^{1/2}
\\
&& \qquad \leqslant c_2n^{-(1+\min\{1/2,\delta/4\})},
\end{eqnarray*}
for $k\geqslant\varepsilon_0n$, for sufficiently large $n$ and for all
$v\in\tau_\mu({\mathbf{M}})$ such that $|v|\leqslant r$, so that
\eqref{eqn5} holds.

To show \eqref{eqn6}, we note that there are positive constants $a,b,c$
independent of $n$ and $k$ such that, for given $V^n_k=v$,
\begin{eqnarray*}
&&\bigl\llVert \bigl(V^n_{k+1}-V^n_k
\bigr) \bigl(V^n_{k+1}-V^n_k
\bigr)^\top\bigr\rrVert
\\
&& \qquad \leqslant \frac{a}{n}\bigl|\exp^{-1}_{\mu}(X_1)\bigr|^2+
\frac{b}{k^2}\bigl(c+\| H_{\mu
,X_1}\|^2\bigr)|v|^2.
\end{eqnarray*}
Thus, by result (i),
\begin{eqnarray*}
&&\bigl\llVert \mathrm{E}\bigl[\bigl(V^n_{k+1}-V^n_k
\bigr) \bigl(V^n_{k+1}-V^n_k
\bigr)^\top 1_{\{ |V^n_{k+1}-V^n_k|>1\}}\vert V^n_k=v\bigr]
\bigr\rrVert
\\
&& \qquad \leqslant  \frac{a}{n}\mathrm{E}\bigl[\bigl|\exp^{-1}_{\mu
}(X_1)\bigr|^21_{\{\llvert V^n_{k+1}-V^n_k\rrvert >1\}}
\vert V^n_k=v\bigr]+\frac{b}{k^2}\mathrm{E}
\bigl[c+\|H_{\mu ,X_1}\|^2\bigr]|v|^2
\\
&& \qquad \leqslant\frac{a}{n}\mathrm{E}\bigl[\bigl|\exp^{-1}_{\mu
}(X_1)\bigr|^{2+\delta}
\bigr]^{2/(2+\delta
)}\mathrm{P}\bigl(\bigl\llvert V^n_{k+1}-V^n_k
\bigr\rrvert >1\vert V^n_k=v\bigr)^{\delta/(2+\delta
)}
\\
&&\qquad \quad{}+\frac{b}{k^2}\mathrm{E}\bigl[c+\|H_{\mu,X_1}\|^2
\bigr]|v|^2
\\
&& \qquad \leqslant \frac{a}{n}\mathrm{E}\bigl[\bigl|\exp^{-1}_{\mu
}(X_1)\bigr|^{2+\delta}
\bigr]^{2/(2+\delta
)}\times\frac{c'}{n^{\delta/(2+\delta)}}+\frac
{b}{k^2}\mathrm{E}
\bigl[c+\|H_{\mu ,X_1}\|^2\bigr]|v|^2
\end{eqnarray*}
for some constant $c'>0$ dependent on $|v|$, so that the required
result follows.
\end{pf}

\begin{corollary*}\label{cor1}
Under the assumptions of Lemma~\ref{lem2}, for any $\varepsilon_0>0$ and
$r>0$, the following limits hold uniformly in $k\geqslant\varepsilon_0n$:
\begin{longlist}[(iii)]
\item[(i)] for any $\varepsilon>0$,
\[
\lim_{n\rightarrow\infty}\sup_{|v|\leqslant r}n\mathrm{P}\bigl(
\bigl\llvert V^n_{k+1}-V^n_k\bigr\rrvert >\varepsilon\vert V^n_{k}=v\bigr)=0;
\]
\item[(ii)]
\[
\lim_{n\rightarrow\infty}\sup_{|v|\leqslant r}\bigl\llvert n
\mathrm{E}\bigl[ \bigl(V^n_{k+1}-V^n_k
\bigr)1_{\{\llvert V^n_{k+1}-V^n_k\rrvert \leqslant1\} }\vert V^n_{k}=v\bigr]\bigr\rrvert
=0;
\]
\item[(iii)]
\begin{eqnarray*}
&&
\lim_{n\rightarrow\infty}\sup_{|v|\leqslant
r}\bigl\llVert n
\mathrm{E}\bigl[\bigl(V^n_{k+1}-V^n_k
\bigr) \bigl(V^n_{k+1}-V^n_k
\bigr)^\top 1_{\{ |V^n_{k+1}-V^n_k|\leqslant1\}}\vert V^n_k=v
\bigr]-A\bigr\rrVert
\\
&& \qquad =0,
\end{eqnarray*}
where $A=\mathrm{E}[H_{\mu,X_1}]^{-1}\Gamma \mathrm{E}[H_{\mu,X_1}]^{-\top}$ and
%
\begin{equation}\label{eqn8}
\Gamma=\operatorname{cov}\bigl(\exp^{-1}_\mu(X_1)
\bigr)=\mathrm{E}\bigl[\exp ^{-1}_\mu(X_1)
\otimes \exp^{-1}_\mu(X_1)\bigr].
\end{equation}
\end{longlist}
\end{corollary*}

\begin{pf}
By \eqref{eqn4}, for any $k\geqslant\varepsilon_0n$,
\begin{eqnarray*}
&&\mathrm{P}\bigl(\bigl\llvert V^n_{k+1}-V^n_k
\bigr\rrvert >\varepsilon\vert V^n_k=v\bigr)
\\
&& \qquad\leqslant \mathrm{E}\bigl[\rho(\mu,X_1)^{2+\delta}\bigr]
\frac
{2^{2+\delta}}{(\varepsilon
'\sqrt{n})^{2+\delta}}+\mathrm{E}\bigl[\|H_{\mu,X_1}\|^2\bigr]
\frac
{(4|v|)^2}{(\varepsilon
'\varepsilon_0n)^2},
\end{eqnarray*}
when $n$ is sufficiently large. Thus (i) holds. Noting that
$\mathrm{E}[V^n_{k+1}\vert V^n_k]=V^n_k$, (ii)~follows from
\eqref{eqn5}. Since
\[
\operatorname{cov}\bigl(V^n_{k+1}\vert V^n_k
\bigr)=\operatorname{cov} \biggl(\frac{1}{\sqrt
{n}} \bigl(\mathrm{E}[H_{\mu,X_1}]
\bigr)^{-1}\bigl(\exp^{-1}_\mu (X_1)
\bigr) \biggr)=\frac{1}{n}A,
\]
(iii) is equivalent to
\[
\lim_{n\rightarrow\infty}\sup_{|v|>r}\bigl\llVert n
\mathrm{E}\bigl[\bigl(V^n_{k+1}-V^n_k
\bigr) \bigl(V^n_{k+1}-V^n_k
\bigr)^\top1_{\{|V^n_{k+1}-V^n_k|>1\}
}\vert V^n_k=v\bigr]
\bigr\rrVert =0,
\]
which follows from \eqref{eqn6}.
\end{pf}

The properties that we have obtained so far on $\{V^n_k\vert k\geqslant
0\}$ enable us to prove the weak convergence of $\{V^n_{[nt]}\vert t\geqslant0\}$ as follows.

\begin{proposition*}\label{prop}
In addition to the assumptions in \eqref{eqn3}, assume that, for some
$\delta>0$, $E[\rho(\mu,X_1)^{2+\delta}]<\infty$. Then the
sequence of
processes $\{V^n_{[nt]}\vert t\geqslant0\}$ converges weakly in
$\mathbb
{D}([0,\infty),\tau_{\mu}({\mathbf{M}}))$, the space of right continuous
functions with left limits on the tangent space of $\mathbf{M}$ at
$\mu$, to
$\{V_t\vert t\geqslant0\}$ as $n\rightarrow\infty$, where $V_t$ is the
solution of the stochastic differential equation
%
\begin{equation}\label{eqn7}
dV_t= \bigl\{ \bigl(\mathrm{E}[H_{\mu,X_1}]
\bigr)^{-1}\Gamma \bigl(\mathrm{E}[H_{\mu ,X_1}]
\bigr)^{-\top} \bigr\} ^{1/2}\, dB_t
\end{equation}
with $V_0=0$, $B_t$ a standard Brownian motion in $\mathbb{R}^d$ and
$\Gamma$ is defined by \eqref{eqn8}.
\end{proposition*}

\begin{pf}
Let $\tilde V^n_k=(\frac{k}{n},V^n_k)$. Then $\{\tilde V^n_k\vert k\geqslant0\}$ is a time-homogeneous Markov chain. For each
$n\geqslant
1$, write $P_n$ for the transition probability distribution associated
with $\{\tilde V^n_k\vert k\geqslant0\}$, that is,
\[
P_n \biggl( \biggl(\frac{l}{n},v \biggr),B \biggr)=
\mathrm{P}\biggl( \biggl(\frac {l+1}{n},V^n_{l+1}
\biggr)\in B\vert V^n_l=v\biggr),
\]
where $B$ is any Borel set in $(0,\infty)\times\tau_\mu({\mathbf{M}})$.

For any $\varepsilon_0>0$, the result of Lemma~\ref{lem1} implies that
$\{
\tilde V^n_{[\varepsilon_0n]}\vert n\geqslant1\}$ is tight. Hence, there is
a subsequence $\{\tilde V^{n_j}_{[\varepsilon_0{n_j}]}\vert j\geqslant1\}$
that converges weakly in $\tau_\mu({\mathbf{M}})$ to a random variable
$\tilde\xi_{\varepsilon_0}=(\varepsilon_0,\xi_{\varepsilon_0})$. Then it follows
from Corollary~7.4.2 in \cite{EK} (pages 355--356) that the results of
the \hyperref[cor1]{Corollary} imply that the sequence of processes $\{\tilde
V^{n_j}_{[n_jt]}\vert t\geqslant\varepsilon_0\}$ converges weakly in
$\mathbb{D}([\varepsilon_0,\infty),[\varepsilon_0,\infty)\times\tau
_\mu(M))$
to a diffusion $\{\tilde V_t\vert t\geqslant\varepsilon_0\}$, where
$\tilde
V_t=(t,V_t)$ with the initial condition that $\tilde V_{\varepsilon_0}$
has the same distribution as $\tilde\xi_{\varepsilon_0}$ and where $V_t$
satisfies the stochastic differential equation \eqref{eqn7}. This
implies (cf. \cite{EK}, page 355) that $\{V^{n_j}_{[n_jt]}\vert t\geqslant\varepsilon_0\}$ converges weakly in $\mathbb{D}([\varepsilon
_0,\infty),\tau_\mu(M))$ to $\{V_t\vert t\geqslant\varepsilon_0\}$, where
$V_{\varepsilon_0}$ has the same distribution as $\xi_{\varepsilon_0}$.

To show the required result, it is now sufficient to show that, for any
subsequence of $\{V^n_{[nt]}\vert t\geqslant0\}$, there is a further
subsequence which converges weakly to $\{V_t\vert t\geqslant0\}$.
Without loss of generality, we may rename the subsequence as $\{
V_{[nt]}^n\vert t\geqslant0\}$ and apply the above to $\varepsilon_0=1/m$.
For each $m\geqslant1$, this gives a subsequence $\{
V^{n_j}_{[n_jt]}\vert t\geqslant0\}$ indexed by $m$, of $\{
V^n_{[nt]}\vert t\geqslant0\}$, which converges weakly on $[1/m,\infty)$
to $\{V_t\vert t\geqslant1/m\}$. Hence, we obtain a sequence indexed by
$m$ of subsequences, and we then take the diagonal subsequence. The
diagonal subsequence converges weakly in $D((0,\infty),\tau_\mu
({\mathbf{M}}))$ to $\{V'_t\vert t>0\}$. However, the result of Lemma~\ref{lem1}
shows that $\mathrm{E}[|V_t'|^2]\rightarrow0$ as $t\rightarrow
0$ and so the
required result follows by noting that $\{V_t'\vert t\geqslant0\}$ must
be equal in law to $\{V_t\vert t\geqslant0\}$.
\end{pf}

\section{The main result}\label{sec3}

We now return to consider the sample Fr\'echet means of $\{X_k\vert k\geqslant1\}$. For this, we denote by $\mu_k$ a sample Fr\'echet mean
of $X_1,\ldots,X_k$ for each $k$, so that $\mu_k$ converges to $\mu$
almost surely (cf. \cite{HZ}). It follows from \eqref{eqn1} that $\mu
_k$ satisfies the condition
%
\begin{equation}\label{eqn9}
\frac{1}{k}\sum_{i=1}^k
\exp_{\mu_k}^{-1}(X_i)=0.
\end{equation}
Thus the origin of the tangent space of $\mathbf{M}$ at $\mu_k$,
$\tau_{\mu
_k}({\mathbf{M}})$, is the sample Euclidean mean of the Euclidean random
variables $\exp^{-1}_\mu(X_i)$, $i=1,\ldots,k$. Nevertheless, although
these relations resemble those for Euclidean means, these conditions
are generally imposed on different tangent spaces, resulting in the
difficulty in obtaining a usable form of the relation between
consecutive sample Fr\'echet means. Moreover, the usual difference
``$\mu
_k-\mu$'' makes no sense here. However, in the context of manifolds,
$\exp^{-1}_\mu(\mu_k)$ plays a similar role to $\mu_k-\mu$ in the
Euclidean case. This leads us to consider, for each $n\geqslant1$, the
re-scaled sequence
%
\begin{equation}\label{eqn12}
W^n_k=\frac{k}{\sqrt n}\exp_{\mu}^{-1}(
\mu_k),\qquad k\geqslant 1.
\end{equation}
It is clear from \eqref{eqn9}, which the sample Fr\'echet means must
satisfy, that $\mu_{k+1}$ cannot generally be expected to be determined
by $\mu_k$ and $X_{k+1}$ alone so that in particular, $\{\mu_k\vert k\geqslant1\}$, and so $\{W^n_k\vert k\geqslant1\}$, is in general
\textit{not} a Markov chain. However, the following result shows that,
for sufficiently large $n$ and $k$, the behaviour of $\{W^n_k\vert k\geqslant1\}$ is close to that of a Markov chain.

\begin{lemma}\label{lem3}
In addition to the assumptions in \eqref{eqn3}, assume that
\[
\lim_{s\rightarrow0}E \Bigl[\sup_{x\in\operatorname{ball}(\mu,s)}\|
\Pi_{x,\mu
}H_{x,X_1}-H_{\mu,X_1}\| \Bigr]=0,
\]
where $\Pi_{x,y}$ denotes the parallel transport from $x$ to $y$ along
the geodesic between the two points. Then, for any $\varepsilon_0>0$,
$r>0$, and $T>0$,
\[
\sup_{\varepsilon_0\leqslant t\leqslant T\wedge\sigma_n^r}\bigl\llvert \bigl(W^n_{[nt]}-V^n_{[nt]}
\bigr)-\bigl(W^n_{[\varepsilon_0n]}-V^n_{[\varepsilon
_0n]}\bigr)
\bigr\rrvert \stackrel{\mathrm{P}} {\longrightarrow}0 \qquad\mbox{as }n\rightarrow
\infty,
\]
where $\{V^n_k\vert k\geqslant0\}$ are the Markov chains defined in the
previous section and $\sigma_n^r=\inf\{t\geqslant\varepsilon_0\vert |W^n_{[nt]}|\geqslant r\mbox{ or }|W^n_{[nt]-1}|\geqslant r\}$.
\end{lemma}

Note that, when $x$ is sufficiently close to $\mu$, the geodesic
between the two points is unique so that the above parallel transport
is well defined.

Note also that $\operatorname{Hess}_x(\frac{1}{2}\rho(x,y)^2)$ is, as a mapping
from $\tau_x(\mathbf{M})\times\tau_x(\mathbf{M})\mapsto\mathbb
R$, smooth with
respect to $x$ if $y\notin\mathcal{C}_x$ and, by \eqref{eqn0}, it is
positive-definite provided $\sqrt{\kappa_1}\rho(x,y)<\pi/2$. Thus the
relationship between $H_{x,y}$ and $\operatorname{Hess}_x(\frac{1}{2}\rho
(x,\break y)^2)$ ensures that all three assumptions required for Lemma~\ref
{lem3} are satisfied if the support for the distribution of $X$ is a
compact subset of the open ball $\operatorname{ball}(\mu,\pi/(2\sqrt{\kappa_1}))$.

\begin{pf*}{Proof of Lemma~\ref{lem3}}
Define, for each given $k$, the random vector field $U_k$ on~$\mathbf
{M}$ by\vspace*{-6pt}
\[
U_k(x)=\sum_{i=1}^k
\exp^{-1}_x(X_i),
\]
for $x\notin\mathcal{C}_{X_1}\cup\mathcal{C}_{X_2}\cup\cdots
\cup
\mathcal{C}_{X_k}$. For each fixed $x$, $\frac{1}{k}U_k(x)$ is the
sample Euclidean mean of random variables $\exp_x^{-1}(X_1),\ldots
,\exp
_x^{-1}(X_k)$. By hypothesis on $X_i$ and the result of \cite{LB},
$U_k(\mu)$ is defined almost surely, and it follows from \eqref{eqn1}
that $E[U_k(\mu)]=0$. Moreover, $\mu_k$ being a sample Fr\'echet mean
of $X_1,\ldots,X_k$ implies that $U_k(\mu_k)=0$ almost surely. Using
these facts and using parallel transport followed by Taylor's expansion,
Kendall and Le \cite{KL} show that
%
\begin{equation}\label{eqn10}
\quad-\sum_{i=1}^kD_{\exp^{-1}_\mu(\mu_k)}
\exp^{-1}_\mu(X_i)=U_k(\mu )+\Delta
_k(\mu_k;X_1,\ldots,X_k) \bigl(
\exp^{-1}_\mu(\mu_k)\bigr),\hspace*{-10pt}
\end{equation}
where the correction operator $\Delta_k$ satisfies the condition that,
for any given $\varepsilon>0$, there exists $s>0$ such that the ball,
$\operatorname{ball}(\mu,s)$, that is centred at $\mu$ and with radius~$s$ is
contained in ${\mathbf{M}}\setminus\mathcal{C}_\mu$ and, for any
$x$ in that ball,
\begin{eqnarray*}
&&\bigl\|\Delta_k(x;X_1,\ldots,X_k)\bigr\|
\\
&& \qquad \leqslant d\sum_{i=1}^k \Bigl\{(1+2
\varepsilon s)\sup_{x'\in\operatorname{ball}(\mu
,s)}\bigl\|\Pi_{x',\mu}D\exp^{-1}_{x'}(X_i)-D
\exp^{-1}_\mu(X_i)\bigr\|
\\
&&\hspace*{107pt}\qquad \quad{}+2\varepsilon \bigl(\bigl|\exp^{-1}_\mu(X_i)\bigr|+s\bigl\|D
\exp^{-1}_\mu (X_i)\bigr\| \bigr) \Bigr\}.
\end{eqnarray*}
Thus, noting by \eqref{eqn2} that
\[
- (D_{\exp^{-1}_\mu(\mu_k)}U_k ) (\mu)= \Biggl(\sum
_{i=1}^kH_{\mu
,X_i} \Biggr) \bigl(
\exp^{-1}_\mu(\mu_k)\bigr),
\]
we can rewrite \eqref{eqn10} as
%
\begin{equation} \label{eqn11}
\Biggl\{\sum_{i=1}^kH_{\mu,X_i}-
\Delta_k(\mu_k;X_1,\ldots ,X_k)
\Biggr\} \bigl(\exp^{-1}_\mu(\mu_k)
\bigr)=U_k(\mu),
\end{equation}
which leads to a link between $W^n_k$ and the rescaled sample Euclidean
mean $U_k(\mu)$.

On the other hand,
\begin{eqnarray*}
&& \Biggl\{\sum_{i=1}^{k+1}H_{\mu,X_i}-
\Delta_{k+1}(\mu _{k+1};X_1,\ldots
,X_{k+1}) \Biggr\}\bigl(\exp^{-1}_\mu(
\mu_k)\bigr)
\\
&& \qquad = \Biggl\{\sum_{i=1}^kH_{\mu,X_i}-
\Delta_k(\mu_k;X_1,\ldots ,X_k)
\Biggr\} \bigl(\exp^{-1}_\mu(\mu_k)
\bigr)+H_{\mu,X_{k+1}}\bigl(\exp^{-1}_\mu(\mu_k)
\bigr)
\\
&&\qquad\quad{}+ \bigl\{\Delta_k(\mu_k;X_1,
\ldots,X_k)-\Delta_{k+1}(\mu _{k+1};X_1,
\ldots,X_{k+1}) \bigr\}\bigl(\exp^{-1}_\mu(
\mu_k)\bigr)
\\
&& \qquad =U_k(\mu)+H_{\mu,X_{k+1}}\bigl(\exp^{-1}_\mu(
\mu_k)\bigr)+R(X_1,\ldots ,X_{k+1}) \bigl(
\exp^{-1}_\mu(\mu_k)\bigr),
\end{eqnarray*}
where
\[
R(X_1,\ldots,X_{k+1})=\Delta_k(
\mu_k;X_1,\ldots,X_k)-\Delta _{k+1}(
\mu _{k+1};X_1,\ldots,X_{k+1}).
\]
This, together with \eqref{eqn11}, gives
\begin{eqnarray*}
&& \Biggl\{\sum_{i=1}^{k+1}H_{\mu,X_i}-
\Delta_{k+1}(\mu _{k+1};X_1,\ldots
,X_{k+1}) \Biggr\} \bigl(\exp^{-1}_\mu(
\mu_{k+1})- \exp^{-1}_\mu (\mu _k) \bigr)
\\
&& \qquad =\exp^{-1}_\mu(X_{k+1})-H_{\mu,X_{k+1}}\bigl(
\exp^{-1}_\mu(\mu _k)\bigr)-R(X_1,
\ldots,X_{k+1}) \bigl(\exp^{-1}_\mu(
\mu_k)\bigr).
\end{eqnarray*}

It then follows from the definition of $W^n_k$ that the difference
$W^n_{k+1}-W^n_k$ can be expressed as
\begin{eqnarray*}
&&W^n_{k+1}-W^n_k
\\
&& \qquad =\frac{k+1}{\sqrt n} \Biggl\{\sum_{i=1}^{k+1}H_{\mu,X_i}-
\Delta _{k+1}(\mu_{k+1};X_1,\ldots,X_{k+1})
\Biggr\}^{-1}
\\
&& \qquad\quad{}\times \bigl\{\exp^{-1}_\mu(X_{k+1})-H_{\mu,X_{k+1}}
\bigl(\exp ^{-1}_\mu (\mu_k)\bigr)\\
&& \qquad \qquad\hspace*{23pt}{}-R(X_1,
\ldots,X_{k+1}) \bigl(\exp^{-1}_\mu(
\mu_k)\bigr) \bigr\}
\\
&&\qquad\quad{}+\frac{1}{\sqrt n}\exp^{-1}_\mu(\mu_k),
\end{eqnarray*}
or equivalently as
\begin{eqnarray*}
&&W^n_{k+1}- \biggl(1+\frac{1}{k}
\biggr)W^n_k
\\
&& \qquad = \Biggl\{\frac{1}{k+1} \Biggl(\sum_{i=1}^{k+1}H_{\mu,X_i}-
\Delta _{k+1}(\mu_{k+1};X_1,\ldots,X_{k+1})
\Biggr) \Biggr\}^{-1}
\\
&&\qquad\quad {}\times \biggl\{\frac{1}{\sqrt n}\exp^{-1}_\mu(X_{k+1})-
\frac
{1}{k}H_{\mu,X_{k+1}}\bigl(W^n_k\bigr)-
\frac{1}{k}R(X_1,\ldots ,X_{k+1})
\bigl(W^n_k\bigr) \biggr\}.
\end{eqnarray*}
However, under the given assumptions, we have
\[
\frac{1}{k}\sum_{i=1}^k
H_{\mu,X_i}\stackrel{\mathrm{a.s.}}{\longrightarrow}E[H_{\mu,X_1}]
\quad\mbox{and}\quad \frac{1}{k}\bigl\|\Delta_k(\mu_k;X_1,
\ldots,X_k)\bigr\|\stackrel{\mathrm {P}}{\longrightarrow}0
\]
(cf. \cite{KL}), so that in particular, $\frac{1}{k}\|R(X_1,\ldots
,X_{k+1})\|\stackrel{\mathrm{P}}{\longrightarrow}0$. Hence, it
follows that
\begin{eqnarray*}
W^n_{k+1}&=&\frac{1}{\sqrt n}\bigl(E[H_{\mu,X_1}]
\bigr)^{-1}\exp^{-1}_\mu (X_{k+1})
\\
&&{}+ \biggl\{\frac{k+1}{k}I-\frac{1}{k}\bigl(E[H_{\mu,X_1}]
\bigr)^{-1}H_{\mu
,X_{k+1}} \biggr\}\bigl(W^n_k
\bigr)+ o\bigl(k^{-1}\bigr) \qquad \mbox{a.s.},
\end{eqnarray*}
where $I$ is the identity operator. This implies that, for $t\geqslant
\varepsilon_0$,
\begin{eqnarray*}
&&\bigl(W^n_{[nt]}-V^n_{[nt]}\bigr)-
\bigl(W^n_{[\varepsilon_0n]}-V^n_{[\varepsilon_0n]}\bigr)
\\
&& \qquad =\bigl(W^n_{[nt]-1}-V^n_{[nt]-1}\bigr)-
\bigl(W^n_{[\varepsilon_0n]}-V^n_{[\varepsilon
_0n]}\bigr)
\\
&& \qquad = \frac{1}{[nt]-1} \bigl\{I-\bigl(E[H_{\mu,X_1}]\bigr)^{-1}H_{\mu
,X_{[nt]}}
\bigr\}\bigl(W^n_{[nt]-1}-V^n_{[nt]-1}
\bigr)\\
&&\qquad\quad{}+o\bigl([nt]^{-1}\bigr)\qquad \mbox{a.s.}
\end{eqnarray*}
so that, for $\varepsilon_0\leqslant t\leqslant T\wedge\sigma_n^r$,
\begin{eqnarray*}
&&\bigl\llvert \bigl(W^n_{[nt]}-V^n_{[nt]}
\bigr)-\bigl(W^n_{[\varepsilon_0n]}-V^n_{[\varepsilon
_0n]}\bigr)
\bigr\rrvert
\\
&& \qquad \leqslant \bigl\llvert \bigl(W^n_{[nt]-1}-V^n_{[nt]-1}
\bigr)-\bigl(W^n_{[\varepsilon
_0n]}-V^n_{[\varepsilon_0n]}\bigr)
\bigr\rrvert +o\bigl([nt]^{-1}\bigr).
\end{eqnarray*}
The required result then follows.
\end{pf*}

We are now in the position to state and prove the main result of the
paper concerning the limiting diffusion associated with the sequences
of the rescaled images $\{W^n_k\vert k\geqslant0\}$, under $\exp
^{-1}_\mu
$, of the Fr\'echet means $\mu_k$ of $X_1,\ldots,X_k$.

\begin{theorem*}\label{thm1}
Under the assumptions of the \hyperref[prop]{Proposition} and Lemma~\ref{lem3}, the
sequence of processes $\{W^n_{[nt]}\vert t\geqslant0\}$ converges weakly
in $\mathbb{D}([0,\infty),\tau_{\mu}({\mathbf{M}}))$ to $\{V_t\vert t\geqslant
0\}$, where $W^n_0=0$; $W^n_k$, $k\geqslant1$, is defined by \eqref{eqn12}, and the $V_t$ are as given in the \hyperref[prop]{Proposition}.
\end{theorem*}

\begin{pf}
By the \hyperref[prop]{Proposition}, we only need to show that, for any $r>0$ and $T>0$,
\[
\sup_{0\leqslant t\leqslant T\wedge\sigma_n^r}\bigl\llvert W^n_{[nt]}-V^n_{[nt]}
\bigr\rrvert \stackrel{\mathrm{P}} {\longrightarrow }0\qquad \mbox{as }n\rightarrow
\infty,
\]
where $\sigma_n^r=\inf\{t\geqslant\varepsilon_0\vert |W^n_{[nt]}|\geqslant
r\mbox{ or }|W^n_{[nt]-1}|\geqslant r\}$.

Since $W^n_0=V^n_0=0$, we have
\begin{eqnarray*}
&&\sup_{0\leqslant t\leqslant T\wedge\sigma_n^r}\bigl\llvert W^n_{[nt]}-V^n_{[nt]}
\bigr\rrvert
\\
&& \qquad =\lim_{\varepsilon_0\downarrow0}\sup_{\varepsilon_0\leqslant
t\leqslant
T\wedge\sigma_n^r}\bigl\llvert
\bigl(W^n_{[nt]}-V^n_{[nt]}\bigr)\bigr
\rrvert
\\
&& \qquad\leqslant \lim_{\varepsilon_0\downarrow0}\Bigl\{\sup_{\varepsilon
_0\leqslant
t\leqslant T\wedge\sigma_n^r}\bigl
\llvert \bigl(W^n_{[nt]}-V^n_{[nt]}\bigr)-
\bigl(W^n_{[\varepsilon_0n]}-V^n_{[\varepsilon
_0n]}\bigr)\bigr
\rrvert
\\
&&\hspace*{133pt}\qquad\quad{}+\bigl\llvert W^n_{[\varepsilon_0n]}\bigr\rrvert +\bigl\llvert
V^n_{[\varepsilon
_0n]}\bigr\rrvert \Bigr\}.
\end{eqnarray*}
Thus, for any $\varepsilon>0$, we have for all sufficiently small
$\varepsilon_0>0$,
\begin{eqnarray*}
&& \mathrm{P}\Bigl(\sup_{0\leqslant t\leqslant T\wedge
\sigma_n^r}\bigl\llvert
\bigl(W^n_{[nt]}-V^n_{[nt]}\bigr)\bigr
\rrvert >6\varepsilon\Bigr)
\\
&& \qquad  \leqslant  \mathrm{P}\Bigl(\sup_{\varepsilon
_0\leqslant t\leqslant T\wedge\sigma_n^r}\bigl\llvert
\bigl(W^n_{[nt]}-V^n_{[nt]}\bigr)-
\bigl(W^n_{[\varepsilon _0n]}-V^n_{[\varepsilon
_0n]}\bigr)\bigr
\rrvert \\
&&\hspace*{87pt}\qquad\qquad{}+\bigl\llvert W^n_{[\varepsilon_0n]}\bigr\rrvert +\bigl\llvert
V^n_{[\varepsilon_0n]}\bigr\rrvert >3\varepsilon\Bigr)
\\
&& \qquad  \leqslant \mathrm{P}\Bigl(\sup_{\varepsilon
_0\leqslant t\leqslant T\wedge\sigma_n^r}\bigl\llvert
\bigl(W^n_{[nt]}-V^n_{[nt]}\bigr)-
\bigl(W^n_{[\varepsilon _0n]}-V^n_{[\varepsilon
_0n]}\bigr)\bigr
\rrvert >\varepsilon\Bigr)
\\
&&\qquad\quad{}+ \mathrm{P}\bigl(\bigl\llvert W^n_{[\varepsilon_0n]}\bigr
\rrvert >\varepsilon \bigr)+\mathrm{P}\bigl(\bigl\llvert V^n_{[\varepsilon_0n]}
\bigr\rrvert >\varepsilon\bigr).
\end{eqnarray*}
The first term on the right tends to zero as $n\rightarrow\infty$ by
Lemma~\ref{lem3}. It follows from the \hyperref[prop]{Proposition} that the distribution
$\nu_{\varepsilon_0}$ of $V_{\varepsilon_0}$ is Gaussian\vspace*{1pt} with mean zero and
covariance matrix $\varepsilon_0^2\mathrm{E}[H_{\mu
,X_1}]^{-1}\Gamma \mathrm{E}[H_{\mu ,X_1}]^{-\top}$, where
$\Gamma$ is given by \eqref{eqn8}. This implies
that the limiting distribution of $V^n_{[\varepsilon_0n]}$ is $\nu
_{\varepsilon_0}$. Moreover, the result of \cite{KL} implies that $\nu
_{\varepsilon_0}$ is also the limiting distribution of $W^n_{[\varepsilon
_0n]}$. Thus, as $n\rightarrow\infty$, both the second and third terms
on the right are bounded above by var$(|V_{\varepsilon_0}|)/\varepsilon^2$,
so that
\[
\lim_{n\rightarrow\infty}\mathrm{P}\Bigl(\sup_{0\leqslant
t\leqslant T\wedge\sigma _n^r}
\bigl\llvert \bigl(W^n_{[nt]}-V^n_{[nt]}
\bigr)\bigr\rrvert >6\varepsilon\Bigr)\leqslant2\frac
{\operatorname{var}(|V_{\varepsilon_0}|)}{\varepsilon^2}.
\]
Since $\lim_{\varepsilon_0\downarrow0}\operatorname{var}(|V_{\varepsilon
_0}|)=0$, the independence of the left-hand side above on $\varepsilon_0$
then gives the required result.
\end{pf}

It is interesting to note the relationship between the result of the
\hyperref[thm1]{Theorem} and the central limit theorem for Fr\'echet means obtained in
\cite{KL}, in comparison with that between the corresponding results
for Euclidean means. It is also interesting to see the difference
between the limiting diffusion obtained here and that obtained in \cite{ADPY}. The latter should shed some light on the difference between the
asymptotic behaviour of the sample Fr\'echet means and that of the
random sequence obtained using the stochastic gradient algorithm
constructed in \cite{ADPY}.






\printaddresses
\end{document}